\documentclass[a4paper,11pt,twoside,reqno]{amsart}

\usepackage[utf8]{inputenc}
\usepackage[plainpages=false,pdfpagelabels=true]{hyperref}
\usepackage{amssymb,amsthm,wasysym}
\usepackage{slashed}

\numberwithin{equation}{section}

\newtheorem{Satz}{Theorem}[section]

\newtheorem{Lem}[Satz]{Lemma}

\newtheorem{Cor}[Satz]{Corollary}
\theoremstyle{definition}
\newtheorem{Dfn}[Satz]{Definition}
\newtheorem{Bem}[Satz]{Remark}

\allowdisplaybreaks[1]

\renewcommand{\epsilon}{\varepsilon}

\newcommand{\R}{\ensuremath{\mathbb{R}}}
\newcommand{\D}{\slashed{D}}
\newcommand{\p}{\slashed{\partial}}
\usepackage[margin=1in]{geometry}

\title{Energy estimates for the supersymmetric nonlinear sigma model and applications}
\author{Volker Branding}
\date{\today}
\address{TU Wien\\
Institut für diskrete Mathematik und Geometrie\\
Wiedner Hauptstraße 8–10, A-1040 Wien}
\email[]{volker@geometrie.tuwien.ac.at}
\begin{document}

\begin{abstract}
We derive gradient and energy estimates for critical points of the full supersymmetric 
sigma model and discuss several applications.
\end{abstract} 

\maketitle

\section{Introduction and Results}
The full nonlinear supersymmetric \(\sigma\)-model is an important model in modern quantum field theory.
In the physical literature \cite{MR1707282}, \cite{Callan:1989nz} it is usually formulated in terms of supergeometry, which includes the use of Grassmann-valued spinors.
However, taking ordinary instead of Grassmann-valued spinors one can investigate the full nonlinear supersymmetric \(\sigma\)-model
as a geometric variational problem. This study was initiated in \cite{MR2176464}, where the notion of \emph{Dirac-harmonic maps}
was introduced. These form a pair of a map between Riemannian manifolds and a vector spinor.
More precisely, the equations for Dirac-harmonic maps couple the harmonic map equation to spinor fields.
As limiting cases both harmonic maps and harmonic spinors can be obtained. 
In the case of a two-dimensional domain Dirac-harmonic maps belong to the class of conformally invariant variational problems yielding a rich structure.

Many important results for Dirac-harmonic maps have already been established.
This includes the regularity of weak solutions \cite{MR2544729}
and an existence result for uncoupled solutions \cite{MR3070562}.
The boundary value problem for Dirac-harmonic maps is studied in \cite{MR3085099}, \cite{MR3044133}.
The heat-flow for Dirac-harmonic maps was studied recently in \cite{regularized}, \cite{regularized-surface}
and \cite{sun}.

However, to analyze the full nonlinear supersymmetric \(\sigma\)-model one has to go beyond the notion of Dirac-harmonic maps.
Considering an additional two-form in the action functional one is led to \emph{magnetic Dirac-harmonic maps} introduced in \cite{MR3305429}.
\emph{Dirac-harmonic maps} to target spaces \emph{with torsion} are analyzed in \cite{torsion}.
Finally, taking into account a curvature term in the action functional one is led to \emph{Dirac-harmonic maps with curvature term},
which were introduced in \cite{MR2370260}. 

In this note we study general properties of the system of partial differential equations 
that arises as critical points of the full nonlinear supersymmetric \(\sigma\)-model.

This article is organized as follows. In Section 2 we recall the mathematical background that we
are using to perform our analysis. In Section 3 we present an \(\epsilon\)-regularity theorem for the domain being a closed surface
and as an application, we prove the removable singularity theorem for Dirac-harmonic maps with curvature term.
In Section 4 we derive gradient estimates and point out several applications.

\section{The full supersymmetric nonlinear sigma model}
Throughout this article, we assume that \((M,h)\) is a Riemannian spin manifold with spinor bundle \(\Sigma M\), 
for more details about spin geometry see the book \cite{MR1031992}. Moreover, let \((N,g)\) be another Riemannian
manifold and let \(\phi\colon M\to N\) be map. Together with the pullback bundle \(\phi^{-1}TN\) 
we can consider the twisted bundle \(\Sigma M\otimes\phi^{-1}TN\). 
The induced connection on this bundle will be denoted by \(\tilde{\nabla}\). 
Sections \(\psi\in\Gamma(\Sigma M\otimes\phi^{-1}TN)\) in this bundle  are called \emph{vector spinors} and the natural operator acting on them 
is the twisted Dirac operator, denoted by \(\D\).  
This is an elliptic, first order operator, which is self-adjoint with respect to the \(L^2\)-norm.
More precisely, the twisted Dirac operator is given by \(\D=e_\alpha\cdot\tilde{\nabla}_{e_\alpha}\), where \(\{e_\alpha\}\) is
an orthonormal basis of \(TM\) and \(\cdot\) denotes Clifford multiplication. 
We are using the Einstein summation convention, that is we sum over repeated indices.
Clifford multiplication is skew-symmetric, namely
\[
\langle\chi,X\cdot\xi\rangle_{\Sigma M}=-\langle X\cdot\chi,\xi\rangle_{\Sigma M}
\]
for all \(\chi,\xi\in\Gamma(\Sigma M)\) and all \(X\in TM\).
Moreover, the twisted Dirac-operator \(\D\) satisfies the following Weitzenböck formula
\begin{equation}
\label{weitzenboeck}
\D^2\psi=-\tilde{\Delta}\psi+\frac{1}{4}R\psi +\frac{1}{2}e_\alpha\cdot e_\beta\cdot R^N(d\phi(e_\alpha),d\phi(e_\beta))\psi.
\end{equation}
Here, \(\tilde{\Delta}\) denotes the connection Laplacian on \(\Sigma M\otimes\phi^{-1}TN\), \(R\) denotes the scalar curvature on \(M\)
and \(R^N\) is the curvature tensor on \(N\).
This formula can be deduced from the general Weitzenböck formula for twisted Dirac operators,
see \cite{MR1031992}, p.\ 164, Theorem 8.17.

We do not present the full energy functional here but rather focus on its critical points.
These satisfy a coupled system of the following form 
\begin{align}
\label{phi-mfd}
\tau(\phi)=&A(\phi)(d\phi,d\phi)+B(\phi)(d\phi,\psi,\psi)+C(\phi)(\psi,\psi,\psi,\psi), \\
\label{psi-mfd}
\D\psi=&E(\phi)(d\phi)\psi+F(\phi)(\psi,\psi)\psi.
\end{align}
Here, \(\tau(\phi)\in\Gamma(\phi^{-1}TN)\) denotes the tension field of the map \(\phi\) and the other terms represent the analytical structure of 
the right hand side. We will always assume that the endomorphisms \(A,B,C,E\) and \(F\) are bounded.

At some points we will assume that the target manifold \(N\) is isometrically embedded in some \(\R^q\)
by the Nash embedding theorem. Then, we have that \(\phi\colon M\to\R^q\) with \(\phi(x)\in N\). The vector spinor \(\psi\)  becomes
a vector of usual spinors \(\psi^1,\psi^2,\ldots,\psi^q\), more precisely \(\psi\in\Gamma(\Sigma M\otimes T\R^q)\).
The condition that \(\psi\) is along the map \(\phi\) is then encoded as
\[
\sum_{i=1}^q\nu^i\psi^i=0\qquad \text{for any normal vector }\nu \text{ at } \phi(x).
\]

The system \eqref{phi-mfd}, \eqref{psi-mfd} then acquires the form
\begin{align}
\label{phi-global}-\Delta\phi=&\tilde{A}(\phi)(d\phi,d\phi)+\tilde{B}(\phi)(d\phi,\psi,\psi)+\tilde{C}(\phi)(\psi,\psi,\psi,\psi),\\
\label{psi-global}\p\psi=&\tilde{E}(\phi)(d\phi)\psi+\tilde{F}(\phi)(\psi,\psi)\psi.
\end{align}
Here \(\p:=e_\alpha\cdot\nabla^{\Sigma M}_{e_\alpha}\) denotes the usual Dirac-operator acting on sections in \(\psi\in\Gamma(\Sigma M\otimes T\R^q)\).

The quantities \(A,B,C,E\) and \(F\) can be extended to the ambient space (denoted by a tilde) and depend only on geometric data.
However, this does not alter the analytic structure of the right hand side of \eqref{phi-mfd}, \eqref{psi-mfd}.

\begin{Bem}
The regularity of the system \eqref{phi-global}, \eqref{psi-global} is already fully understood.
By now, there are powerful tools available to ensure the smoothness of a system like \eqref{phi-global}, \eqref{psi-global},
see \cite{MR3020100}, \cite{MR2661574} and \cite{MR3333092}. 
However, it should be noted that in order to apply the main result from \cite{MR3020100}
we need a certain antisymmetry of the endomorphism \(A\).
It is quite remarkable that the actual \(A\) from the nonlinear supersymmetric sigma model
has the necessary antisymmetry.
\end{Bem}

\begin{Bem}
In the physical literature the energy functional for the full supersymmetric nonlinear sigma model
is fixed by the requirements of superconformal invariance (conformal invariance + supersymmetry)
and invariance under diffeomorphisms on the domain.
\end{Bem}

\section{Energy estimates and applications}
Throughout this section we assume that the domain \(M\) is a closed Riemannian spin surface.
\subsection{Epsilon Regularity Theorem}
We derive an \(\epsilon\)-regularity Theorem for smooth solutions of the system
\eqref{phi-global}, \eqref{psi-global}.
To this end, we combine the methods for Dirac-harmonic maps from \cite{MR2176464}, Theorem 3.2
and nonlinear Dirac equations from \cite{MR2390834}, Theorem 2.1. To establish the \(\epsilon\)-regularity theorem
we make use of the invariance under scaling of the system \eqref{phi-global}, \eqref{psi-global}.

However, we should not assume that the energy is small globally.
\begin{Lem}
Assume that the pair \((\phi,\psi)\) is a smooth solution of	
\eqref{phi-global} and \eqref{psi-global} satisfying
\begin{equation}
\int_M(|d\phi|^2+|\psi|^4)<\epsilon_0
\end{equation}
with \(\epsilon_0\) small enough. Moreover, assume that there are no harmonic spinors on \(M\).
Then both \(\phi\) and \(\psi\) are trivial.
\end{Lem}
\begin{proof}
See the proof of Lem. 4.8. in \cite{MR3333092}.
\end{proof}

We define the following local energy:
\begin{Dfn}
Let \(U\) be a domain on \(M\). We define the energy of the pair \((\phi,\psi)\) on \(U\) by
\begin{equation}
E(\phi,\psi,U):=\int_U(|d\phi|^2+|\psi|^4).
\end{equation}
\end{Dfn}

Similar as in the case of Dirac-harmonic maps \cite{MR2176464} we prove the following
\begin{Satz}[\(\epsilon\)-Regularity Theorem]
\label{epsilon-regularity}
Assume that the pair \((\phi,\psi)\) is a smooth solution of
\eqref{phi-global} and \eqref{psi-global} with small energy 
\begin{equation}
\label{small-energy}
E(\phi,\psi,D)<\epsilon.
\end{equation}
Then the following estimate holds
\begin{align}
|d\phi|_{W^{1,p}(\tilde{D})}&\leq C(\tilde{D},p)(|d\phi|_{L^2(D)}+|\psi|^2_{L^4(D)}),  \\
|\nabla\psi|_{W^{1,p}(\tilde{D})}&\leq C(\tilde{D},p)|\psi|_{L^4(D)}
\end{align}
for all \(\tilde{D}\subset D, p>1\), where \(C(\tilde{D},p)\) is a positive
constant depending only on \(\tilde{D}\) and \(p\).
\end{Satz}

We divide the proof into several steps, we will assume that
\(\tilde{D}\subset D^3\subset D^2\subset D^1\subset D\).

As a first step, we derive an estimate for the spinor \(\psi\), similar to Lemma 3.4 in \cite{MR2176464}.
\begin{Lem}
Assume that the pair \((\phi,\psi)\) is a smooth solution of
\eqref{phi-global} and \eqref{psi-global} satisfying \eqref{small-energy}. 
Then the following estimate holds
\begin{align}
\label{psi-lq}
|\psi|_{L^q(D^1)}&\leq C(D^1)|\psi|_{L^4(D)},\qquad \forall q>1,\qquad \forall D^1\subset D
\end{align}
where \(C(D^1)\) is a constant depending only on \(D^1\).
\end{Lem}

\begin{proof}
We choose a cut-off function \(\eta\) satisfying \(0\leq\eta\leq 1\) with \(\eta|_{D^1}=1\) and \(\operatorname{supp}\eta\subset D\).
Consider the spinor \(\xi:=\eta\psi\) and moreover, since the unit disc \(D\) is flat, we have \(\p^2=-\Delta\).
Using \(\eqref{psi-global}\), we calculate
\begin{align}
\label{xi-lq}
\p(\eta\psi)=\eta\p\psi+\nabla\eta\cdot\psi 
=\eta \tilde{E}(\phi)(d\phi)\psi+\eta \tilde{F}(\phi)(\psi,\psi)\psi+\nabla\eta\cdot\psi.
\end{align}
Hence, employing elliptic estimates we get
\[
|\xi|_{W^{1,q}(D)}\leq C(\big||d\phi||\eta\psi|\big|_{L^q(E)}+\big|\eta|\psi|^3\big|_{L^q(F)}+|\psi|_{L^q(D)}). 
\]
By Hölder's inequality we can estimate
\begin{align*}
\big||d\phi||\eta\psi|\big|_{L^q(D)}\leq&|d\phi|_{L^2(D)}|\eta\psi|_{L^{q^*}(D)},\\
\big|\eta|\psi|^3\big|_{L^q(D)}\leq&|\psi|^2_{L^4(D)}|\eta\psi|_{L^{q^*}(D)}
\end{align*}
with the conjugate Sobolev index \(q^*=\frac{2q}{2-q}\).
By the Sobolev embedding theorem we may then follow
\[
|\xi|_{L^{q^*}(D)}\leq C(\sqrt{E(\phi,\psi)}|\xi|_{L^{q^*}(D)}+|\psi|_{L^q(D)}). 
\]
Thus, if the energy \(E(\phi,\psi)\) is small enough, we have
\[
|\xi|_{L^{q^*}(D)}\leq C|\psi|_{L^q(D)}.
\]
At this point for any \(p>1\) one can always find some \(q<2\) such that \(p=q^*\)
and this yields the first claim.
\end{proof}

\begin{Lem}
Assume that the pair \((\phi,\psi)\) is a smooth solution of
\eqref{phi-global} and \eqref{psi-global} satisfying \eqref{small-energy}. 
Then the following estimate holds
\begin{equation}
\label{phi-w14}
|\phi|_{W^{1,4}(D^2)}\leq C(D^2)\sqrt{\epsilon},\qquad \forall D^2\text{ with } D^2\subset D,
\end{equation}
where the constant \(C\) depends only on \(D^2\).
\end{Lem}
\begin{proof}
Suppose that \(D^2\subset D\).
We choose a cut-off function \(\eta:0\leq\eta\leq 1\) with \(\eta|_{D^2}=1\) and \(\operatorname{supp}\eta\subset D\).
By equation (\ref{phi-global}) we have
\begin{align*}
|\Delta(\eta\phi)|\leq&C(|\phi|+|d\phi|+|d\phi||d(\eta\phi)|+|\phi d\eta|+\big|\eta|d\phi||\psi|^2\big|+\big|\eta|\psi|^4\big|)\\
\leq& C(|\phi|+|d\phi|+|d\phi||d(\eta\phi)|+\big|\eta|d\phi||\psi|^2\big|+\big|\eta|\psi|^4\big|).
\end{align*}
Hence, for any \(p>1\) we have
\begin{equation}
\label{laplace-u-lp}
|\Delta(\eta\phi)|_{L^p}\leq C(\big||d\phi||d(\eta\phi)|\big|_{L^p}+|d\phi|_{L^{p}}
+\big|\eta|d\phi||\psi|^2\big|_{L^p}+\big|\eta|\psi|^4\big|_{L^p}).
\end{equation}
Choosing \(p=\frac{4}{3}\) on the disc \(D\), we find
\[
|\Delta(\eta\phi)|_{L^{\frac{4}{3}}(D)}\leq C(\big||d\phi||d(\eta\phi)|\big|_{L^{\frac{4}{3}}(D)}+|d\phi|_{L^{\frac{4}{3}}(D)}
+\big|\eta|d\phi||\psi|^2\big|_{L^{\frac{4}{3}}(D)}+\big|\eta|\psi|^4\big|_{L^{\frac{4}{3}}(D)}).
\]
Without loss of generality we assume \(\int_D\phi=0\) such that \(|\phi|_{W^{1,p}(D)}\leq C|d\phi|_{L^p(D)}\) for any \(p>0\).
Moreover, by Hölder's inequality we have 
\[
\big||d\phi||d(\eta\phi)|\big|_{L^{\frac{4}{3}}(D)}\leq |\eta\phi|_{W^{1,4}(D)}|d\phi|_{L^2(D)}
\]
such that we may conclude 
\[
|\eta\phi|_{W^{2,\frac{4}{3}}(D)}\leq C(|\eta\phi|_{W^{1,4}(D)}|d\phi|_{L^2(D)}+|d\phi|_{L^{\frac{4}{3}}(D)}
+\big|\eta|d\phi||\psi|^2\big|_{L^{\frac{4}{3}}(D)}+\big|\eta|\psi|^4\big|_{L^{\frac{4}{3}}(D)}).
\]
By the Sobolev embedding theorem we find \(|\eta\phi|_{W^{1,4}(D)}\leq c|\eta\phi|_{W^{2,\frac{4}{3}}(D)}\) and we may follow
\begin{equation}
\label{estimate-etaphi}
(c^{-1}-C|d\phi|_{L^2(D)})|\eta\phi|_{W^{1,4}(D)}\leq C(|d\phi|_{L^\frac{4}{3}(D)}
+\big|\eta|d\phi||\psi|^2\big|_{L^{\frac{4}{3}}(D)}+\big|\eta|\psi|^4\big|_{L^{\frac{4}{3}}(D)}).
\end{equation}
Regarding the last two terms in \eqref{estimate-etaphi} we note that using \eqref{psi-lq}
\begin{align*}
\big|\eta|d\phi||\psi|^2\big|_{L^{\frac{4}{3}}(D)} \leq & C(|\psi|^2_{L^4(D)}|\eta\phi|_{W^{1,4}(D)}+|\psi|^2_{L^4(D)}), \\
\big|\eta|\psi|^4\big|_{L^{\frac{4}{3}}(D)}\leq &\big|\eta|\psi|^2\big|_{L^4(D)}|\psi|^2_{L^4(D)}\leq C|\psi|^2_{L^4(D)}.
\end{align*}
Applying these estimates and choosing \(\epsilon\) small enough, (\ref{estimate-etaphi}) gives
\[
|\eta\phi|_{W^{1,4}(D)}\leq C( |d\phi|_{L^{\frac{4}{3}}(D)}+\sqrt{\epsilon}|\eta\phi|_{W^{1,4}(D)}+|\psi|^2_{L^4(D)}),
\]
which can be rearranged as
\[
|\eta\phi|_{W^{1,4}(D)}\leq C(|d\phi|_{L^\frac{4}{3}(D)}+|\psi|^2_{L^4(D)})\leq \sqrt{\epsilon} C.
\]
Finally, by the properties of \(\eta\) we have that for some \(\epsilon>0\)
\begin{equation}
\label{estimate-step1}
|\phi|_{W^{1,4}(D^2)}\leq C(D^2)(|d\phi|_{L^2(D)}+|\psi|^2_{L^4(D)})\leq \sqrt{\epsilon} C,\qquad \forall D^2\subset D
\end{equation}
holds.
\end{proof}

\begin{Lem}
Assume that the pair \((\phi,\psi)\) is a smooth solution of
\eqref{phi-global} and \eqref{psi-global} satisfying \eqref{small-energy}. 
Then the following estimate holds
\begin{align}
\label{nablapsi-l2}
|\nabla\psi|_{L^2(D^2)}&\leq C(D^2)|\psi|_{L^4(D)},\qquad \forall D^2\subset D,
\end{align}
where \(C(D^2)\) is a constant depending only on \(D^2\).
\end{Lem}

\begin{proof}
We choose a cut-off function \(\eta\) satisfying \(0\leq\eta\leq 1\) with \(\eta|_{D^2}=1\) and \(\operatorname{supp}\eta\subset D\).
Again, consider the spinor \(\xi:=\eta\psi\) and using \eqref{xi-lq} we estimate
\begin{align*}
|\nabla\xi|_{L^2(D)}&\leq C(|\eta\psi|^3_{L^6(D)}+\big|\eta|d\phi||\psi|\big|_{L^2(D)}+|\psi|_{L^2(D)}) \\
&\leq C(|\psi|^3_{L^4(D)}+|\eta d\phi|_{L^4(D)}|\psi|_{L^4(D)}+|\psi|_{L^4(D)}) \\
&\leq C|\psi|_{L^4(D)}(1+|\psi|^2_{L^4(D)}+|d\phi|_{L^4(D^2)}) \\
&\leq C|\psi|_{L^4(D)},
\end{align*}
which proves the claim.
\end{proof}

\begin{Lem}
Assume that the pair \((\phi,\psi)\) is a smooth solution of
\eqref{phi-global} and \eqref{psi-global} satisfying \eqref{small-energy}.  
Then the following estimate holds
\begin{align}
\label{dphi-l4}
|d\phi|_{L^4(D^2)}&\leq C(D^2)(|d\phi|_{L^2(D)}+|\psi|^2_{L^4(D)}),\qquad \forall D^2\text{ with } D^2\subset D
\end{align}
where \(C\) is a constant depending only on \(D^2\).
\end{Lem}

\begin{proof}
Choose a cut-off function \(\eta\colon 0\leq\eta\leq 1\) with \(\eta|_{D^2}=1\) and \(\operatorname{supp}\eta\subset D\).
By \eqref{estimate-etaphi} we have	
\begin{align*}
|\eta\phi|_{W^{1,4}(D)}&\leq C(|d\phi|_{L^\frac{4}{3}(D)}+\big|\eta|\psi|^2|d\phi|\big|_{L^\frac{4}{3}(D)}+\big|\eta|\psi|^4\big|_{L^{\frac{4}{3}}(D)}).
\end{align*}
Using
\begin{align*}
\big||\psi|^2|d\phi|\big|_{L^\frac{4}{3}(D)}&\leq|\psi|^2_{L^8(D)}|d\phi|_{L^2(D)}\leq C|\psi|_{L^4(D)}|d\phi|_{L^2(D)},  \\
\big||\psi|^4\big|_{L^{\frac{4}{3}}(D)}&\leq |\psi|^2_{L^8(D)}|\psi|^2_{L^4(D)}\leq C|\psi|_{L^4(D)}|\psi|^2_{L^4(D)}
\end{align*}
we obtain the result.
\end{proof}

\begin{Lem}
Assume that the pair \((\phi,\psi)\) is a smooth solution of
\eqref{phi-global} and \eqref{psi-global} satisfying \eqref{small-energy}. 
Then the following estimate holds
\begin{equation}
\label{phi-w2p}
|\phi|_{W^{2,p}(D^3)}\leq C(|d\phi|_{L^2(D)}+|\psi|^2_{L^4(D)}),\qquad \forall D^3\subset D,
\end{equation}
where the constant \(C\) depends only on \(D^3\).
\end{Lem}
\begin{proof}
Choose a cut-off function \(\eta\colon 0\leq\eta\leq 1\) with \(\eta|_{D^3}=1\) and \(\operatorname{supp}\eta\subset D^2\).
By (\ref{laplace-u-lp}) we have
\begin{align*}
|\eta\phi|_{W^{2,2}(D^2)} \leq & C(|d(\eta\phi)|_{L^4(D^2)}|d\phi|_{L^4(D^2)}+|\phi|_{W^{1,2}(D^2)}+\big||d\phi||\psi|^2\big|_{L^2(D^2)}
+\big||\psi|^4\big|_{L^2(D^2)}) \\
\leq &C(|\eta\phi|_{W^{1,4}(D^2)}|d\phi|_{L^4(D^2)}+|\phi|_{W^{1,2}(D^2)} 
+|d\phi|_{L^4(D^2)}|\psi|^2_{L^8(D^2)}+|\psi|^4_{L^8(D^2)}).
\end{align*}
By the Sobolev embedding theorem we get
\[
|\eta\phi|_{W^{1,4}(D^2)}\leq c|\eta\phi|_{W^{2,\frac{4}{3}}(D^2)}\leq c|\eta\phi|_{W^{2,2}(D^2)}.
\]
Moreover, applying
\[
|\eta\phi|_{W^{1,4}(D^2)}|d\phi|_{L^4(D^2)}\leq c\sqrt{\epsilon}|\eta\phi|_{W^{2,2}(D^2)}
\]
we find
\begin{align*}
(1-c\sqrt{\epsilon})|\eta\phi|_{W^{2,2}(D^2)} &\leq C(|\phi|_{W^{1,2}(D^2)}+|d\phi|_{L^4(D^2)}|\psi|^2_{L^8(D^2)}+|\psi|^4_{L^8(D^2)}) \\
&\leq C(|\phi|_{W^{1,4}(D^2)}+|\psi|^4_{L^8(D^2)}).
\end{align*}
Hence, we may conclude
\[
|\eta\phi|_{W^{2,2}(D^2)}\leq C(|\phi|_{W^{1,4}(D^2)}+|\psi|^4_{L^8(D^2)})\leq C(|d\phi|_{L^4(D^2)}+|\psi|^2_{L^4(D^2)}).
\]
Again, by the Sobolev embedding theorem we may thus follow
\begin{equation}
\label{phi-lp}
|d\phi|_{L^p(D^3)}\leq C(|d\phi|_{L^4(D^2)}+|\psi|^2_{L^4(D^2)}),\qquad \forall p>1.
\end{equation}
Having gained control over the \(W^{2,2}\) norm of \(\phi\) we now may control the \(W^{2,p}\) norm of \(\phi\) for \(p>2\).
Again, suppose that \(\tilde{D}\subset D^3\) and choose a cut-off function \(\eta\colon 0\leq\eta\leq 1\)
with \(\eta|_{\tilde{D}}=1\) and \(\operatorname{supp}\eta\subset\tilde{D}\). By (\ref{laplace-u-lp}) we have for any \(p>1\)
\[
|\eta\phi|_{W^{2,p}(D^3)}\leq C(\big||d\phi||d(\eta\phi)|\big|_{L^p(D^3)}+|\phi|_{W^{1,p}(D^3)}+\big||d\phi||\psi|^2\big|_{L^p(D^3)}
+\big||\psi|^4\big|_{L^p(D^3)}).
\]
By application of (\ref{phi-lp}) we find
\begin{align*}
\big||d\phi||d(\eta\phi)|\big|_{L^p(D^3)}&\leq |d\phi|^2_{L^{2p}(D^3)}\leq C(|d\phi|_{L^4(D^2)}+|\psi|^2_{L^4(D^2)}), \\
\big||d\phi||\psi|^2\big|_{L^p(D^3)}&\leq |\psi|^2_{L^{4p}(D^3)}|d\phi|_{L^{2p}(D^3)}\leq C(|d\phi|_{L^4(D^2)}+|\psi|^2_{L^4(D^2)}), \\
\big||\psi|^4\big|_{L^p(D^3)}&=|\psi|^4_{L^{4p}(D^3)}\leq C|\psi|^2_{L^4(D^2)},
\end{align*}
which gives
\[
|\eta\phi|_{W^{2,p}(D^3)}\leq C(|d\phi|_{L^4(D^2)}+|\psi|^2_{L^4(D^2)}).
\]
Finally, we conclude by (\ref{dphi-l4}) that
\[
|\phi|_{W^{2,p}(\tilde{D})}\leq C(|d\phi|_{L^4(D^2)}+|\psi|^2_{L^4(D^2)})\leq C(|d\phi|_{L^2(D)}+|\psi|^2_{L^4(D)}),
\]
which proves the assertion.
\end{proof}

After having gained control over \(\phi\) we may now control the spinor \(\psi\).

\begin{Lem}
Assume that the pair \((\phi,\psi)\) is a smooth solution of
\eqref{phi-global} and \eqref{psi-global} satisfying \eqref{small-energy}. 
Then the following estimates hold:
\begin{align}
|\psi|_{L^\infty(D^2)}& \leq C|\psi|_{L^4(D)}, \qquad \forall D^2\text{ with } D^2\subset D,\\
|\nabla\psi|_{L^\infty(D^2)} & \leq C|\psi|_{L^4(D)}, \qquad \forall D^2\text{ with } D^2\subset D,
\end{align}
where the constants depend only on \(D^2\).
\end{Lem}

\begin{proof}
First of all, we calculate
\begin{align*}
-\Delta\psi=\p^2\psi=\p\big(\tilde{E}(\phi)(d\phi)\psi+\tilde{F}(\phi)(\psi,\psi)\psi\big).
\end{align*}
By a direct calculation this leads to
\begin{align*}
\p(\tilde{E}(\phi)(d\phi)\psi)
=&e_\alpha\cdot(\nabla_{d\phi(e_\alpha)}\tilde{E}(\phi))(d\phi)\psi
+e_\alpha\cdot \tilde{E}(\phi)(\nabla_{e_\alpha} d\phi)\psi+\tilde{E}(\phi)(d\phi)\p\psi
\end{align*}
and in addition
\[
\p(\tilde{F}(\phi)(\psi,\psi)\psi)=e_\alpha\cdot(\nabla_{d\phi(e_\alpha)}\tilde{F}(\phi))(\psi,\psi)\psi
+2e_\alpha\cdot \tilde{F}(\phi)(\nabla_{e_\alpha}\psi,\psi)\psi+e_\alpha\cdot\tilde{F}(\phi)(\psi,\psi)\nabla_{e_\alpha}\psi.
\]
Consequently, for any \(\eta\in C^\infty(D,\R)\) with \(0\leq\eta\leq 1\), we may follow
\[
|\Delta(\eta\psi)|\leq C(|\psi|+|\nabla\psi|+|d\phi|^2|\psi|+|d\phi||\nabla\psi|+|\nabla^2\phi||\psi|+|d\phi||\psi|^3+|\nabla\psi||\psi|^2).
\]
Now for \(D^2\subset D^1\) choose a cut-off function \(\eta\colon 0\leq\eta\leq 1\) with \(\eta|_{D^2}=1\) and \(\operatorname{supp}\eta\subset D^1\).
For any \(p>1\) we then have
\begin{align}
\label{psi-w2p}
\big|\eta\psi\big|_{W^{2,p}(D^1)} \leq C&\big(|\psi|_{L^p(D^1)}+|\nabla\psi|_{L^p(D^1)}+\big||d\phi|^2|\psi|\big|_{L^p(D^1)}+\big||d\phi||\nabla\psi|\big|_{L^p(D^1)} \\
\nonumber&+\big||\nabla^2\phi||\psi|\big|_{L^p(D^1)} +\big||d\phi||\psi|^3\big|_{L^p(D^1)}+\big||\nabla\psi||\psi|^2\big|_{L^p(D^1)}
\big).
\end{align}
Setting \(p=\frac{4}{3}\) and making using of Hölder's inequality we obtain
\begin{align*}
\big|\eta\psi\big|_{W^{2,\frac{4}{3}}(D^1)} \leq& C\big(|\psi|_{L^\frac{4}{3}(D^1)}+|\nabla\psi|_{L^\frac{4}{3}(D^1)}+|d\phi|_{L^4(D^1)}^2|\psi|_{L^4(D^1)}+|d\phi|_{L^4(D^1)}|\nabla\psi|_{L^2(D^1)}\\
&+|\nabla^2\phi|_{L^2(D^1)}|\psi|_{L^4(D^1)}+|d\phi|_{L^4(D^1)}|\psi|^3_{L^6(D^1)} +|\nabla\psi|_{L^2(D^1)}|\psi|^2_{L^8(D^1)}	 
\big).
\end{align*}
By application of (\ref{psi-lq}), (\ref{nablapsi-l2}), (\ref{dphi-l4}) and (\ref{phi-w2p}) we get
\[
|\psi|_{W^{2,\frac{4}{3}}(D^2)}\leq C|\psi|_{L^4(D)}.
\]
By the Sobolev embedding theorem this yields
\begin{equation}
\big|\psi|_{W^{1,4}(D^2)}\leq C|\psi|_{L^4(D)} 
\end{equation}
and also
\[
|\psi|_{L^\infty(D^2)}\leq C|\psi|_{L^4(D)}.
\]
This proves the first estimate for the spinor.
\par\medskip
Using the same method as before, we now get an estimate on \(|\nabla\psi|\).
Thus, for \(D^3\subset D^2\) choose a cut-off function \(\eta\colon 0\leq\eta\leq 1\) with \(\eta|_{D^3}=1\) and \(\operatorname{supp}\eta\subset D^2\).
Setting \(p=2\) in (\ref{psi-w2p}) we obtain
\begin{align*}
\big|\eta\psi\big|_{W^{2,2}(D^2)} \leq& C\big(|\psi|_{L^2(D^2)}+|\nabla\psi|_{L^2(D^2)}+\big||d\phi|^2|\psi|\big|_{L^2(D^2)}+\big||d\phi||\nabla\psi|\big|_{L^2(D^2)} \\
&+\big||\nabla^2\phi||\psi|\big|_{L^2(D^2)}+\big||d\phi||\psi|^3\big|_{L^2(D^2)}+\big||\nabla\psi||\psi|^2\big|_{L^2(D^2)}
\big) \\
\leq &  C \big(|\psi|_{L^4(D^2)}+|\psi|_{L^4(D)}+|d\phi|^2_{L^8(D^2)}|\psi|_{L^4(D^2)}+|d\phi|_{L^4(D^2)}|\nabla\psi|_{L^4(D^2)}\\
&+|\nabla^2\phi|_{L^4(D^2)}|\psi|_{L^4(D^2)}+|d\phi|_{L^4(D^2)}|\psi|^3_{L^{12}(D^2)}+|\psi|_{L^4(D)}\big) \\
\leq&  C|\psi|_{L^4(D)}.
\end{align*}
By the Sobolev embedding theorem we may then follow
\begin{equation}
\label{psi-w1p}
|\psi|_{W^{1,p}(D^3)}\leq C|\psi|_{L^4(D^2)}.
\end{equation}
At this point for \(\tilde{D}\subset D^3\) we again use (\ref{psi-w2p}) with a cut-off function \(\eta\colon 0\leq\eta\leq 1\) with \(\eta|_{\tilde D}=1\) and \(\operatorname{supp}\subset D^3\).
Using (\ref{psi-w1p}), (\ref{psi-lq}), (\ref{dphi-l4}) and (\ref{phi-w2p}) we can follow
\[
|\eta\psi|_{W^{2,p}(D^3)}\leq C|\psi|_{L^4(D)},\qquad\forall p>1.
\]
Thus
\[
|\nabla\psi|_{W^{1,p}(\tilde{D})}\leq C|\psi|_{L^4(D)}
\]
and, finally, we obtain
\[
|\nabla\psi|_{L^\infty(\tilde{D})}\leq C|\psi|_{L^4(D)}.
\]
\end{proof}
This proves Theorem \ref{epsilon-regularity}.
By scaling we obtain the following (similar to Cor. 4.4 in \cite{MR2262709})
\begin{Cor}
\label{corollary-energy-local}
There is an \(\epsilon>0\) small enough such that if the pair \((\phi,\psi)\) is a smooth solution of (\ref{phi-global})
and (\ref{psi-global}) on \(D\setminus\{0\}\) with finite energy \(E(\phi,\psi,D)<\epsilon\), 
then for any \(x\in D_{\frac{1}{2}}\) we have
\begin{align}
|d\phi(x)||x|\leq& C(|d\phi|_{L^2(D_{2|x|})}+|\psi|_{L^4(D_{2|x|})}), \\
|\psi(x)|^\frac{1}{2}|x|^\frac{1}{2}+|\nabla\psi(x)||x|^\frac{3}{2}\leq& C|\psi|_{L^4(D_{2|x|})}.
\end{align}
\end{Cor}
\begin{proof}
This follows from a scaling argument, fix any \(x_0\in D\setminus\{0\}\) and define \((\tilde{\phi},\tilde{\psi})\) by
\[
\tilde{\phi}(x):=\phi(x_0+|x_0|x) \textrm{ and } \tilde{\psi}(x):=|x_0|^\frac{1}{2}\psi(x_0+|x_0|x).
\]
It is easy to see that \((\tilde{\phi},\tilde{\psi})\) is a smooth solution of (\ref{phi-global}) and (\ref{psi-global}) on \(D\)
with \(E(\tilde{\phi},\tilde{\psi},D)<\epsilon\). By application of Theorem \ref{epsilon-regularity}, we have
\begin{align*}
|d\tilde{\phi}|_{L^\infty(D_\frac{1}{2})}\leq C(|d\tilde{\phi}|_{L^2(D)}+|\tilde{\psi}|_{L^4(D)}), \qquad
|\tilde{\psi}|_{C^1(D_\frac{1}{2})}\leq C|\tilde{\psi}|_{L^4(D)}
\end{align*}
and scaling back yields the assertion.
\end{proof}

\subsection{Application: Removable Singularity Theorem for Dirac-harmonic maps with curvature term}
Using the previous estimates we sketch how to prove the removable singularity theorem for Dirac-harmonic maps
with curvature term.

Dirac-harmonic maps with curvature term are critical points of the functional
\begin{equation}
\label{energy-functional}
E_c(\phi,\psi)=\frac{1}{2}\int_M|d\phi|^2+\langle\psi,\D\psi\rangle-\frac{1}{6}\langle R^N(\psi,\psi)\psi,\psi\rangle
\end{equation}
with the indices contracted as 
\[
\langle R^N(\psi,\psi)\psi,\psi\rangle=R_{ijkl}\langle\psi^i,\psi^k\rangle\langle\psi^j,\psi^l\rangle.
\]
The critical points of the energy functional \eqref{energy-functional} are given by (see \cite{MR3333092}, Prop. 2.1)
\begin{align}
\label{dhc-phi}\tau(\phi)=&\frac{1}{2}R^N(e_\alpha\cdot\psi,\psi)d\phi(e_\alpha)-\frac{1}{12}\langle(\nabla R^N)^\sharp(\psi,\psi)\psi,\psi\rangle, \\
\label{dhc-psi}\D\psi=&\frac{1}{3}R^N(\psi,\psi)\psi,
\end{align}
where \(\tau(\phi)\) is the tension field of the map \(\phi\), \(R^N\) denotes the curvature tensor on \(N\)
and \(\sharp\colon\phi^{-1}T^\ast N\to\phi^{-1}TN\) represents the musical isomorphism.

By embedding \(N\) into \(\R^q\) isometrically the equations \eqref{dhc-phi} and \eqref{dhc-psi} acquire the form \eqref{phi-global} and \eqref{psi-global}.
For more details see Lemma 3.5 in \cite{MR3333092}.
\begin{Lem}
\label{lemma-polar-coordinates}
Let \((\phi,\psi)\) be a smooth Dirac-harmonic map with curvature term on \(D\setminus\{0\}\)
satisfying \(E(\phi,\psi,D)<\epsilon\). Then we have
\begin{align}
\int_0^{2\pi}\frac{1}{r^2}|\phi_\theta|^2d\theta&=\int_0^{2\pi}(|\phi_r|^2+\langle\psi,\partial_r\cdot\frac{\tilde{\nabla}\psi}{\partial r}\rangle
-\frac{1}{3}(1+\sin^2\theta)\langle R^N(\psi,\psi)\psi,\psi\rangle)d\theta\\
\nonumber &=\int_0^{2\pi}(|\phi_r|^2-\frac{1}{r^2}\langle\psi,\partial_\theta\cdot\frac{\tilde{\nabla}\psi}{\partial\theta}\rangle
-\frac{\sin^2\theta}{3}\langle R^N(\psi,\psi)\psi,\psi\rangle)d\theta,
\end{align}
where \((r,\theta)\) are polar coordinates on the disc \(D\) around the origin, \(\phi_r\) denotes differentiation of \(\phi\)
with respect to \(r\) and \(\phi_\theta\) denotes differentiation of \(\phi\)
with respect to \(\theta\).
\end{Lem}
 
\begin{proof}
On a small domain \(\tilde{M}\) of \(M\) we choose a local isothermal parameter \(z=x+iy\) and set
\begin{align}
\label{hopf-differential}
T(z)dz^2=&(|\phi_x|^2-|\phi_y|^2-2i\langle\phi_x,\phi_y\rangle 
+\langle\psi,\partial_x\cdot\tilde{\nabla}_{\partial_x}\psi\rangle-i\langle\psi,\partial_x\cdot\tilde{\nabla}_{\partial_y}\psi\rangle 
-\frac{1}{3}\langle R^N(\psi,\psi)\psi,\psi\rangle)dz^2
\end{align}
with \(\partial_x=\frac{\partial}{\partial x}\) and \(\partial_y=\frac{\partial}{\partial y}\).
It was shown in \cite{MR3333092}, Prop. 3.3, that the quadratic differential \eqref{hopf-differential} is holomorphic.
By Corollary \ref{corollary-energy-local} we know that
\begin{align*}
|d\phi|^2\leq\frac{C}{z^2},\qquad
|\psi||\tilde{\nabla}\psi|\leq C(|\psi||\nabla\psi|+|d\phi||\psi|^2)\leq\frac{C}{z^2},\qquad
|\langle R^N(\psi,\psi)\psi,\psi\rangle|\leq\frac{C}{z^2},
\end{align*}
which, altogether gives \(|T(z)|\leq Cz^{-2}\).
Moreover, it is easy to see that \(\int_D|T(z)|<\infty\).
Hence, we may follow that \(zT(z)\) is holomorphic on the disc \(D\) and by Cauchy's integral theorem we deduce
\begin{equation}
0=\operatorname{Im}\int_{|z|=r}zT(z)dz=\int_0^{2\pi}\operatorname{Re}(z^2T(z))d\theta.
\end{equation}
By a direct calculation we find
\begin{align*}
\langle\psi,\partial_x\cdot\tilde{\nabla}_{\partial_x}\psi\rangle-i\langle\psi,\partial_x\cdot\tilde{\nabla}_{\partial_y}\psi\rangle=&
\cos^2\theta\langle\psi,\partial_r\cdot\tilde{\nabla}_{\partial_r}\psi\rangle
-\frac{\sin^2\theta}{r^2}\langle\psi,\partial_\theta\cdot\tilde{\nabla}_{\partial_\theta}\psi\rangle \\
&+\frac{\sin\theta\cos\theta}{r}(\langle\psi,\partial_r\cdot\tilde{\nabla}_{\partial_\theta}\psi\rangle-\langle\psi,\partial_\theta\cdot\tilde{\nabla}_{\partial_r}\psi\rangle).
\end{align*} 
Using the equation for \(\psi\) in polar coordinates 
\begin{equation}
\label{psi-polarcoordinates}
\partial_r\cdot\tilde{\nabla}_{\partial_r}\psi+\frac{1}{r^2}\partial_\theta\cdot\tilde{\nabla}_{\partial_\theta}\psi=\frac{1}{3}R^N(\psi,\psi)\psi 
\end{equation}
we find that the term
\[
\langle\psi,\partial_r\cdot\tilde{\nabla}_{\partial_\theta}\psi\rangle-\langle\psi,\partial_\theta\cdot\tilde{\nabla}_{\partial_r}\psi\rangle
=\frac{r^2}{3}\langle\psi,\partial_r\cdot\partial_\theta\cdot R^N(\psi,\psi)\psi\rangle
\]
is both purely real and purely imaginary and thus vanishes.
Thus, we obtain
\begin{align*}
\operatorname{Re}(z^2T(z))=&r^2|\phi_r|^2-|\phi_\theta|^2
+r^2\cos^2\theta\langle\psi,\partial_r\cdot\tilde{\nabla}_{\partial_r}\psi\rangle
-\sin^2\theta\langle\psi,\partial_\theta\cdot\tilde{\nabla}_{\partial_\theta}\psi\rangle \\
&-\frac{r^2}{3}\langle R^N(\psi,\psi)\psi,\psi\rangle,
\end{align*}
which together with \eqref{psi-polarcoordinates} proves the result.
\end{proof}
 
\begin{Satz}[Removable Singularity Theorem]
Let \((\phi,\psi)\) be a Dirac-harmonic map with curvature term which is smooth on 
\(U\setminus\{p\}\) for some \(p\in U\subset M\). If the pair \((\phi,\psi)\)
has finite energy, then \((\phi,\psi)\) extends to a smooth solution on \(U\).
\end{Satz}
\begin{proof}
We do not give a full proof here.
Using the \(\epsilon\)-regularity Theorem \ref{epsilon-regularity} and Lemma \ref{lemma-polar-coordinates}
the removable singularity theorem can be proven the same way as for Dirac-harmonic maps, see the
proof of Theorem 4.6 in \cite{MR2262709} and the proof of Theorem 3.1 in \cite{MR2390834}.
\end{proof}

\section{Gradient estimates and applications}
In this section we derive gradient estimates for solutions \((\phi,\psi)\) of the coupled system \eqref{phi-mfd}, \eqref{psi-mfd}.
To achieve this we extend the techniques from \cite{MR3144359} and \cite{MR647905}, see also \cite{MR2370260}.

\begin{Bem}
In this section we do not necessarily have to assume that the domain \(M\) is compact.
Moreover, we do not have to restrict to a two-dimensional domain \(M\).
However, in the case of the nonlinear supersymmetric sigma model
the term \(A(d\phi,d\phi)\) originates from the variation of a two-form.
If we would assume that \(m=\dim M\geq 2\) then this term would be proportional to \(|d\phi|^m\).
\end{Bem}

To derive a gradient estimate for solutions of \eqref{phi-mfd} and \eqref{psi-mfd}, we recall the following Bochner formula 
for a map \(\phi\colon M\to N\), that is
\[
\Delta\frac{1}{2}|d\phi|^2=|\nabla d\phi|^2+\langle d\phi(\text{Ric}^M(e_\alpha)),d\phi(e_\alpha)\rangle
-\langle R^N(d\phi(e_\alpha),d\phi(e_\beta))d\phi(e_\alpha),d\phi(e_\beta)\rangle
+\langle\nabla\tau(\phi),d\phi\rangle.
\]
Using \eqref{phi-mfd} and by a direct calculation we find
\begin{align*}
\langle\nabla\tau(\phi),d\phi\rangle=&\langle(\nabla_{d\phi}A(\phi))(d\phi,d\phi),d\phi\rangle+2\langle A(\phi)(\nabla d\phi,d\phi),d\phi\rangle \\
\nonumber &+\langle(\nabla_{d\phi}B(\phi))(d\phi,\psi,\psi),d\phi\rangle+\langle B(\phi)(\nabla d\phi,\psi,\psi),d\phi\rangle
+2\langle B(\phi)(d\phi,\tilde{\nabla}\psi,\psi),d\phi\rangle \\
\nonumber &+\langle(\nabla_{d\phi}C(\phi))(\psi,\psi,\psi,\psi),d\phi\rangle+4\langle C(\phi)(\tilde{\nabla}\psi,\psi,\psi,\psi),d\phi\rangle
\end{align*}
and thus we may estimate
\begin{align*}
\Delta\frac{1}{2}|d\phi|^2\geq&|\nabla d\phi|^2 -\kappa_1|d\phi|^2-\kappa_2|d\phi|^4-c_1|d\phi|^4-2c_2|\nabla d\phi||d\phi|^2 -c_3|d\phi|^3|\psi|^2 \\
\nonumber &-c_4|\nabla d\phi||\psi|^2|d\phi|-2c_4|d\phi|^2|\psi||\tilde{\nabla}\psi|-c_5|\psi|^4|d\phi|^2-4c_6|\tilde{\nabla}\psi||\psi|^3|d\phi|
\end{align*}
with the constants \(\text{Ric}^M\geq-\kappa_1, K^N\leq \kappa_2, c_1:=|\nabla A|_{L^\infty}, c_2:=|A|_{L^\infty}, c_3:=|\nabla B|_{L^\infty}, \\
c_4:=|B|_{L^\infty}, c_5:=|\nabla C|_{L^\infty}, c_6:=|C|_{L^\infty}\). Here, \(K^N\) denotes the sectional curvature on \(N\).
Hence, we may rearrange
\begin{align}
\label{dphi-first}
\Delta\frac{1}{2}|d\phi|^2\geq & (1-\delta_2-\delta_4)|\nabla d\phi|^2 -(\kappa_2+c_1+\frac{c_2^2}{\delta_2}+\delta_3+\frac{c_4^2}{\delta_4})|d\phi|^4-\kappa_1|d\phi|^2 \\
\nonumber &-(\frac{c^2_3}{4\delta_3}+\frac{c_4^2}{4\delta_4}+c_5+\frac{4c_6^2}{\delta_6})|d\phi|^2|\psi|^4-(\delta_4+\delta_6)|\psi|^2|\tilde{\nabla}\psi|^2,
\end{align}
where \(\delta_i,i=2,3,4,6\) are positive constants to the determined later.
As a next step we derive an estimate for \(\Delta|\psi|^4\).
By a direct calculation we obtain (with \(R\) being the scalar curvature on \(M\))
\begin{align*}
\Delta\frac{1}{2}|\psi|^4=&2|\psi|^2|\tilde{\nabla}\psi|^2+\big|d|\psi|^2\big|+\frac{R}{2}|\psi|^4+|\psi|^2\langle e_\alpha\cdot e_\beta\cdot R^N(d\phi(e_\alpha),d\phi(e_\beta))\psi,\psi\rangle \\
&-2|\psi|^2\langle\psi,\D^2\psi\rangle,
\end{align*}
where we applied \eqref{weitzenboeck}. To estimate the last term, we use the equation for \(\psi\), \eqref{psi-mfd}, and find
\begin{align*}
\langle\psi,\D (E(\phi)(d\phi)\psi)\rangle=&\langle\psi,\tilde{\nabla}(E(\phi)(d\phi))\cdot\psi\rangle+\langle\psi,E(\phi)(d\phi)\D\psi\rangle, \\
\langle\psi,\D (F(\phi)(\psi,\psi)\psi)\rangle=&\langle\psi,\tilde{\nabla}(F(\phi)(\psi,\psi)\psi))\cdot\psi\rangle+\langle\psi,e_\alpha\cdot F(\phi)(\psi,\psi)\tilde{\nabla}_{e_\alpha}\psi\rangle.
\end{align*}
Due to the skew-symmetry of the Clifford multiplication the first terms on the right hand side are both purely imaginary and purely real and thus vanish.
Moreover, we have the estimate 
\begin{align*}
-2|\psi|^2|\langle\psi,\D^2\psi\rangle|\geq& -2|E|_{L^\infty}|\psi|^3|\D\psi||d\phi|-2|F|_{L^\infty}|\psi|^5\sqrt{m}|\tilde{\nabla}\psi|  \\
\geq&-2\sqrt{m}|E|_{L^\infty}|\psi|^3|\tilde{\nabla}\psi||d\phi|-2\sqrt{m}|F|_{L^\infty}|\psi|^5|\tilde{\nabla}\psi|.
\end{align*}
Again, we may rearrange
\begin{align}
\label{psi-first}
\Delta\frac{1}{2}|\psi|^4\geq& \big|d|\psi|^2\big|^2+\frac{R}{2}|\psi|^4+(2-\delta_7-\delta_8)|\psi|^2|\tilde{\nabla}\psi|^2 
-(m\kappa_3+\frac{mc_7^2}{\delta_7})|\psi|^4|d\phi|^2-\frac{mc_8^2}{\delta_8}|\psi|^8
\end{align}
with the constants \(\kappa_3:=|R^N|_{L^\infty}, c_7:=|E|_{L^\infty}\) and \(c_8:=|F|_{L^\infty}\).
Moreover, \(\delta_7\) and \(\delta_8\) are positive constants to be determined later.
We set
\begin{equation}
e(\phi,\psi):=\frac{1}{2}(|d\phi|^2+|\psi|^4)
\end{equation}
and in addition \(t:=\delta_2+\delta_4\).
Adding up \eqref{dphi-first} and \eqref{psi-first} we obtain
\begin{align}
\label{laplace-energy}
\Delta e(\phi,\psi)\geq &
(1-t)|\nabla d\phi|^2+\big|d|\psi|^2\big|^2+(2-\delta_4-\delta_6-\delta_7-\delta_8)|\psi|^2|\tilde{\nabla}\psi|^2 \\
\nonumber &-(\kappa_2+c_1+\frac{c_2^2}{\delta_2}+\delta_3+\frac{c_4^2}{\delta_4})|d\phi|^4-\kappa_1|d\phi|^2\\
\nonumber &-(\frac{c^2_3}{4\delta_3}+\frac{c_4^2}{4\delta_4}+c_5+\frac{4c_6^2}{\delta_6}+m\kappa_3+\frac{mc_7^2}{\delta_7})|d\phi|^2|\psi|^4 
+\frac{R}{2}|\psi|^4-\frac{mc_8^2}{\delta_8}|\psi|^8.
\end{align}
This allows us to derive a first (similar to \cite{MR631097} for harmonic maps and \cite{MR3144359} for Dirac-harmonic maps)
\begin{Satz}
Let \((\phi,\psi)\) be a smooth solution of \eqref{phi-mfd} and \eqref{psi-mfd}. 
Suppose that \(M\) is a closed Riemannian manifold with positive Ricci curvature
and that the sectional curvature of \(N\) is bounded. If
\begin{equation}
e(\phi,\psi)<\epsilon 
\end{equation}
for \(\epsilon\) small enough, then \(\phi\) is constant and \(\psi\) vanishes identically.
\end{Satz}
\begin{proof}
We use \eqref{laplace-energy}, set \(\delta_4+\delta_6+\delta_7+\delta_8=2\) and \(t=1\).
Then we obtain the estimate
\begin{align*}
\Delta e(\phi,\psi)\geq \big(\kappa_1-\tilde{c}_1|d\phi|^2-\tilde{c}_2|\psi|^4\big)|d\phi|^2+\big(\frac{R}{2}-\frac{mc_8^2}{\delta_8}\big)|\psi|^4,
\end{align*}
where \(\tilde{c}_1>0\) and \(\tilde{c}_2>0\) can be determined from \eqref{laplace-energy} and the above choices for the \(\delta_i,i=2,4,6,7,8\).
By assumption the domain \(M\) has positive Ricci curvature, thus \(\kappa_1\) and \(R\) are both positive. Hence, for \(\epsilon\)
small enough the energy \(e(\phi,\psi)\) is a subharmonic function, which proves the result.
\end{proof}
For the sake of completeness we give the following
\begin{Lem}
We have the following inequality:
\begin{equation}
\label{kato-inequality}
\frac{|de(\phi,\psi)|^2}{2e(\phi,\psi)}\leq |\nabla d\phi|^2+\big|d|\psi|^2\big|^2.
\end{equation}
\end{Lem}
\begin{proof}
We follow \cite{MR3144359}, p.73.
We calculate
\[
de(\phi,\psi)=\langle d\phi,\nabla d\phi\rangle+|\psi|^2d|\psi|^2 
\]
and by squaring the equation we obtain
\begin{align*}
|de(\phi,\psi)|^2\leq &|d\phi|^2|\nabla d\phi|^2+|\psi|^4\big|d|\psi|^2\big|^2+2|d\phi||\nabla d\phi||\psi|^2\big|d|\psi|^2\big| \\
\leq &(|d\phi|^2+|\psi|^4)|\nabla d\phi|^2+(|d\phi|^2+|\psi|^4)\big|d|\psi|^2\big|^2\\
=&2e(|\nabla d\phi|^2+\big|d|\psi|^2\big|^2)
\end{align*}
yielding the result.
\end{proof}

\begin{Lem}
Let \((\phi,\psi)\) be a smooth solution of \eqref{phi-mfd} and \eqref{psi-mfd}. 
Moreover, suppose that the Ricci-curvature of \(M\) satisfies \(\operatorname{Ric}^M\geq-\kappa_1\)
and the sectional curvature \(K^N\) of \(N\) satisfies \(K^N\leq\kappa_2\).
Then the following inequality holds:
\begin{equation}
\label{energy-inequality}
\frac{\Delta e(\phi,\psi)}{e(\phi,\psi)}\geq\frac{1-t}{2}\frac{|de(\phi,\psi)|^2}{e(\phi,\psi)^2}-\frac{m}{2}\kappa_1-c_{13}|d\phi|^2-c_{14}|\psi|^4,
\end{equation}
where the value of the positive constants \(c_{13}\) and \(c_{14}\) is determined along the proof.
\end{Lem}
\begin{proof}
We choose \(\delta_j,j=2,4,6,7,8\) such that
\[
2-\delta_4-\delta_6-\delta_7-\delta_8>0
\]
and \(1-t>0\).
Using \eqref{kato-inequality} we find
\begin{align}
\label{energy-inequality-1}
\Delta e(\phi,\psi)\geq & \frac{1-t}{2}\frac{|de(\phi,\psi)|^2}{e(\phi,\psi)}-\kappa_1|d\phi|^2+\frac{R}{2}|\psi|^4 
-c_{10}|d\phi|^4-c_{11}|d\phi|^2|\psi|^4 -c_{12}|\psi|^8
\end{align}
with the positive constants
\begin{align*}
c_{10}:=&\kappa_2+c_1+\frac{c_2^2}{\delta_2}+\delta_3+\frac{c_4^2}{\delta_4},\qquad
c_{11}:=\frac{c^2_3}{4\delta_3}+\frac{c_4^2}{4\delta_4}+c_5+\frac{4c_6^2}{\delta_6}+m\kappa_3+\frac{mc_7^2}{\delta_7}, \qquad
c_{12}:=\frac{mc_8^2}{\delta_8}.
\end{align*}
Since \(\operatorname{Ric}\geq -\kappa_1\) we have \(R\geq-m\kappa_1\). Dividing by \(e(\phi,\psi)\), using that
\begin{align*}
-2c_{10}\frac{|d\phi|^4}{|d\phi|^2+|\psi|^4}&> -2c_{10}|d\phi|^2,\qquad -2c_{11}\frac{|d\phi|^2|\psi|^4}{|d\phi|^2+|\psi|^4}> -2c_{11}|\psi|^4, \\
-2c_{12}\frac{|\psi|^8}{|d\phi|^2+|\psi|^4}&> -2c_{12}|\psi|^4
\end{align*}
and setting \(c_{13}:=2c_{10}, c_{14}:=2c_{11}+2c_{12}\) we obtain the result.
\end{proof}

\begin{Bem}
If we set 
\[
C:=\min(c_{10},\frac{c_{11}}{2},c_{12})
\]
in \eqref{energy-inequality-1} then we would get an inequality of the form
\[
\frac{\Delta e(\phi,\psi)}{e(\phi,\psi)}\geq\frac{1-t}{2}\frac{|de(\phi,\psi)|^2}{e(\phi,\psi)^2}-\frac{m}{2}\kappa_1-Ce(\phi,\psi).
\]
This energy inequality has the same analytic structure as in the case of harmonic maps.
\end{Bem}

To obtain a gradient estimate from \eqref{energy-inequality} for non-compact \(M\) and \(N\) we need the following tools:
Let \(\rho\) be the Riemannian distance function from the point \(y_0\) in the target manifold \(N\). We define 
\begin{equation}
\xi:=\sqrt{d_1}\cos(\sqrt{d_1}\rho) 
\end{equation}
for some positive number \(\sqrt{d_1}\) to be fixed later, where \(B_R(y_0)\) denotes the geodesic ball
of radius \(R\) around the point \(y_0\).
We will assume that \(R<\pi/(2\sqrt{d_1})\), thus \(0<\xi(R)<\sqrt{d_1}\) on the ball \(B_R(y_0)\). 
\begin{Lem}
On the geodesic ball \(B_R(y_0)\) we have the following estimate
\begin{equation}
\operatorname{Hess}\xi\leq-d_1^\frac{3}{2}\cos(\sqrt{d_1}\rho).
\end{equation}
\end{Lem}
\begin{proof}
This follows from  the Hessian Comparison theorem, see \cite{MR521983}, p.19, Prop. 2.20 
and \cite{MR647905}, p.93.
\end{proof}

In addition, let \(r\) be the distance function from the point \(x_0\) in \(M\). 
Define the function
\begin{equation}
\label{definition-F}
F:=\frac{a^2-r^2}{\xi\circ\phi}e(\phi,\psi)^{p}
\end{equation}
on the geodesic ball \(B_r(x_0)\) in \(M\) with some positive number \(p\).
Clearly, the function \(F\) vanishes on the boundary \(B_a(x_0)\),
hence \(F\) attains its maximum at an interior point \(x_{max}\).
Moreover, we can assume that the distance function \(r\) is smooth near the point \(x_{max}\),
see \cite{MR573431}, section 2.

\begin{Lem}
Suppose that \((M,h)\) and \((N,g)\) are complete Riemannian manifolds. 
Let \((\phi,\psi)\) be a smooth solution of \eqref{phi-mfd} and \eqref{psi-mfd}
satisfying \(\phi\colon M\to B_{R}(y_0)\subset N\) with \(R<\pi/(2\sqrt{d_1})\).
Moreover, suppose that the Ricci-curvature of \(M\) satisfies \(\operatorname{Ric}^M\geq-\kappa_1\)
and the sectional curvature \(K^N\) of \(N\) satisfies \(K^N\leq\kappa_2\).
Then the following inequality holds:
\begin{align}
0\geq&\frac{-\Delta r^2}{a^2-r^2}-\big(1+\frac{1+t}{2}\frac{1}{p}\big)\frac{|d(r^2)|^2}{(a^2-r^2)^2}
-p\frac{m}{2}\kappa_1-pc_{13}|d\phi|^2-pc_{14}|\psi|^4 \\
&\nonumber-(1+t)\frac{|d(r^2)||d(\xi\circ\phi)|}{p(a^2-r^2)\xi\circ\phi}
+\big(1-\frac{1+t}{2}\frac{1}{p}\big)\frac{|d(\xi\circ\phi)|^2}{(\xi\circ\phi)^2}
-\frac{\Delta(\xi\circ\phi)}{\xi\circ\phi}
\end{align}
\end{Lem}
\begin{proof}
Differentiating \(\log F\) at its maximum \(x_{max}\) we obtain
\begin{equation}
\label{F-maximum}
0=\frac{-d(r^2)}{a^2-r^2}-\frac{d(\xi\circ\phi)}{\xi\circ\phi}+p\frac{de(\phi,\psi)}{e(\phi,\psi)}
\end{equation}
and also
\begin{equation}
\label{F-maximum-laplace}
0\geq\frac{-\Delta r^2}{a^2-r^2}-\frac{|d(r^2)|^2}{(a^2-r^2)^2}
-\frac{\Delta(\xi\circ\phi)}{\xi\circ\phi}+\frac{|d(\xi\circ\phi)|^2}{(\xi\circ\phi)^2}
+p\frac{\Delta e(\phi,\psi)}{e(\phi,\psi)}-p\frac{|de(\phi,\psi)|^2}{e(\phi,\psi)^2}.
\end{equation}
Inserting \eqref{energy-inequality} into \eqref{F-maximum-laplace} we find
\begin{align}
\label{F-1}
0\geq&\frac{-\Delta r^2}{a^2-r^2}-\frac{|d(r^2)|^2}{(a^2-r^2)^2}-p\frac{m}{2}\kappa_1-pc_{13}|d\phi|^2-pc_{14}|\psi|^4-p\frac{1+t}{2}\frac{|de(\phi,\psi)|^2}{e(\phi,\psi)^2} \\
\nonumber &-\frac{\Delta(\xi\circ\phi)}{\xi\circ\phi}+\frac{|d(\xi\circ\phi)|^2}{(\xi\circ\phi)^2}.
\end{align}
By squaring \eqref{F-maximum} we also get
\begin{equation}
\label{F-maximum-squared}
p\frac{|d(e(\phi,\psi))|^2}{e(\phi,\psi)^2}\leq\frac{1}{p}\frac{|d(r^2)|^2}{(a^2-r^2)^2}+\frac{2|d(r^2)||d(\xi\circ\phi)|}{p(a^2-r^2)\xi\circ\phi}+\frac{1}{p}\frac{|d(\xi\circ\phi)|^2}{(\xi\circ\phi)^2}.
\end{equation}
Combining \eqref{F-maximum-squared} and \eqref{F-1} then gives the result.
\end{proof}

In the following, we apply the Laplacian comparison Theorem, see \cite{MR521983}, p.20, that is
\[
\Delta r^2\leq C_L(1+r)
\]
with some positive constant \(C_L\). Moreover, we make use of the Gauss Lemma, that is \(|dr|^2=1\).

\begin{Cor}
Suppose that \((M,h)\) and \((N,g)\) are complete Riemannian manifolds. 
Let \((\phi,\psi)\) be a smooth solution of \eqref{phi-mfd} and \eqref{psi-mfd}
satisfying \(\phi\colon M\to B_{R}(y_0)\subset N\) with \(R<\pi/(2\sqrt{d_1})\).
Moreover, suppose that the Ricci-curvature of \(M\) satisfies \(\operatorname{Ric}^M\geq-\kappa_1\)
and the sectional curvature \(K^N\) of \(N\) satisfies \(K^N\leq\kappa_2\).
Then the following inequality holds:
\begin{align}
\label{energy-inequality-cor}
0\geq&\frac{-C_L(1+r)}{a^2-r^2}-\big(1+\frac{1+t}{2}\frac{1}{p}\big)\frac{4r^2}{(a^2-r^2)^2}
-p\frac{m}{2}\kappa_1-pc_{13}|d\phi|^2-pc_{14}|\psi|^4 \\
&\nonumber-(1+t)\frac{2r|d(\xi\circ\phi)|}{p(a^2-r^2)\xi\circ\phi}
+\big(1-\frac{1+t}{2}\frac{1}{p}\big)\frac{|d(\xi\circ\phi)|^2}{(\xi\circ\phi)^2}
-\frac{\operatorname{Hess}\xi(d\phi,d\phi)}{\xi\circ\phi}-\frac{d\xi(\tau(\phi))}{\xi\circ\phi}
\end{align}
\end{Cor}
\begin{proof}
This follows from the Laplacian comparison Theorem, the Gauss Lemma and the chain rule for the tension field
of composite maps, that is
\[
\Delta(\xi\circ\phi)=\operatorname{Hess}\xi(d\phi,d\phi)+d\xi(\tau(\phi)).
\]
\end{proof}

To shorten the notation, we set
\begin{equation}
\label{definition-L_1}
L_1:=\frac{C_L(1+r)}{a^2-r^2}+\big(1+\frac{1+t}{2}\frac{1}{p}\big)\frac{4r^2}{(a^2-r^2)^2}+p\frac{m}{2}\kappa_1.
\end{equation}
By assumption the map \(\phi\) satisfies the equation \eqref{phi-mfd}.
Hence, we may estimate
\[
|\tau(\phi)|\leq|A|_{L^\infty}|d\phi|^2+|B|_{L^\infty}|d\phi||\psi|^2+|C|_{L^\infty}|\psi|^4\leq c_2|d\phi|^2+c_4|d\phi||\psi|^2+c_6|\psi|^4.
\]
Moreover, we have \(|d\xi|=d_1|\sin(\sqrt{d_1}\rho)|\leq d_1\) and
to obtain a gradient estimate we set
\begin{align*}
p=\frac{1+t}{2}=\frac{1+\delta_2+\delta_4}{2}.
\end{align*}
By the properties of the Riemannian distance function \(\rho\) on \(N\),
equation \eqref{energy-inequality-cor}, the definition of \(L_1\) and the estimate on \(\operatorname{Hess}\xi\) we find
\begin{align}
\label{estimate-a}
0\geq &-L_1-\frac{4rd_1}{(a^2-r^2)\xi\circ\phi}|d\phi|+\big(d_1-(1+\delta_2+\delta_4)(\kappa_2+c_1+\frac{c_2^2}{\delta_2}+\delta_3+\frac{c_4^2}{\delta_4})\big)|d\phi|^2 \\
\nonumber&-\frac{1+\delta_2+\delta_4}{2}c_{14}|\psi|^4-\frac{d_1}{\xi\circ\phi}(c_2|d\phi|^2+c_4|d\phi||\psi|^2+c_6|\psi|^4).
\end{align}
\begin{Bem}
If we consider the limiting case of harmonic maps in \eqref{estimate-a} then we obtain the same 
inequality leading to a gradient estimate as in \cite{MR647905}.
\end{Bem}
First of all, let us consider the case that \(A(d\phi,d\phi)=0\) in \eqref{phi-mfd}, which means that \(c_1=c_2=\delta_1=\delta_2=0\) and we obtain
\begin{align}
\label{estimate-b}
0\geq &-L_1-\frac{4rd_1}{(a^2-r^2)\xi\circ\phi}|d\phi|+\big(d_1-(1+\delta_4)(\kappa_2+\delta_3+\frac{c_4^2}{\delta_4})-\delta_9\big)|d\phi|^2 \\
\nonumber&-\big(\frac{1+\delta_4}{2}c_{14}+\frac{d_1^2c_4^2}{4(\xi\circ\phi)^2\delta_9}+\frac{d_1c_6}{\xi\circ\phi}\big)|\psi|^4
\end{align}
for some positive number \(\delta_9\).
We require the coefficient in front of \(|d\phi|^2\) to be positive, which in this case can be expressed as
\begin{align}
\label{definition-d}
\tilde{d}:=d_1-(1+\delta_4)(\kappa_2+\delta_3+\frac{c_4^2}{\delta_4})-\delta_9>0.
\end{align}
Hence, we have to choose \(d_1\) such that \eqref{definition-d} holds.
However, note that we have some freedom to choose \(\delta_3\) and \(\delta_9\) in \eqref{definition-d}.
Again, to shorten the notation, we set
\begin{equation}
\label{definition-L_2}
L_2:=\frac{1+\delta_4}{2}c_{14}+\frac{d_1^2c_4^2}{4(\xi\circ\phi)^2\delta_9}+\frac{d_1c_6}{\xi\circ\phi}
\end{equation}
and then \eqref{estimate-b} becomes
\begin{align}
\label{estimate-c}
0>\tilde{d}|d\phi|^2-\frac{4rd_1}{(a^2-r^2)\xi\circ\phi}|d\phi|-L_1-L_2|\psi|^4.
\end{align}
Note that \eqref{estimate-c} is an equation for an unknown \(x\) of the form
\[
0>ax^2-bx-c
\]
with the constants \(a,b,c\) being all positive. Then, it follows directly that
\[
x<\frac{b}{a}+\sqrt{\frac{c}{a}},
\]
which gives us the following
\begin{Satz}
Suppose that \((M,h)\) and \((N,g)\) are complete Riemannian manifolds. 
Let \((\phi,\psi)\) be a smooth solution of \eqref{phi-mfd} and \eqref{psi-mfd} satisfying
\(\phi\colon M\to B_{R}(y_0)\subset N\) with \(R<\pi/(2\sqrt{d_1})\), where \(d_1\) is determined by \eqref{definition-d}.
Suppose that \(A(d\phi,d\phi)=0\) and \(B,C,E,F\) are bounded. Moreover, assume that the Ricci curvature of \(M\) satisfies
\(\operatorname{Ric}\geq-\kappa_1\) and that the sectional curvature \(K^N\) of \(N\) satisfies \(K^N\leq\kappa_2\).
Then for any \(x_0\in B_a(x_0)\) the following estimate holds
\begin{equation}
\label{gradient-estimate-dphi}
|d\phi|\leq\frac{4rd_1}{\tilde{d}(a^2-r^2)\xi\circ\phi}+\sqrt{\frac{L_1+L_2|\psi|^4}{\tilde{d}}},
\end{equation}
where \(L_1\) is given by \eqref{definition-L_1} and \(L_2\) is given by \eqref{definition-L_2}.
\end{Satz}
In the case that \(A(d\phi,d\phi)\neq 0\) it is more difficult to obtain an estimate on \(|d\phi|\). Let us again consider \eqref{estimate-a}
\begin{align}
\label{estimate-d}
0\geq &-L_1-\frac{4rd_1}{(a^2-r^2)\xi\circ\phi}|d\phi|+\big(d_1-(1+\delta_2+\delta_4)(\kappa_2+c_1+\frac{c_2^2}{\delta_2}+\delta_3+\frac{c_4^2}{\delta_4})-\delta_{10}-\frac{c_2d_1}{\xi\circ\phi}\big)|d\phi|^2 \\
\nonumber&-\big(\frac{1+\delta_4+\delta_2}{2}c_{14}+\frac{d_1^2c_4^2}{4(\xi\circ\phi)^2\delta_{10}}+\frac{d_1c_6}{\xi\circ\phi}\big)|\psi|^4
\end{align}
for some positive number \(\delta_{10}\).
Again, we require the coefficient in front of \(|d\phi|^2\) to be positive, which in this case can be expressed as
\begin{equation}
\label{d-general}
\tilde{d}:=d_1-(1+\delta_2+\delta_4)(\kappa_2+c_1+\frac{c_2^2}{\delta_2}+\delta_3+\frac{c_4^2}{\delta_4})-\delta_{10}-\frac{c_2\sqrt{d_1}}{\cos(\sqrt{d_1}R)}>0.
\end{equation}
However, it seems quite difficult to check if one can arrange all the constants above such that the inequality \eqref{d-general} holds.

\begin{Bem}
\begin{enumerate}
\item Due to the additional terms on the right hand side of \eqref{psi-mfd} it is hard to say in which cases the estimate \eqref{gradient-estimate-dphi} 
is sharp.
\item It becomes clear along the proof that we have a lot of freedom rearranging the constants involved in
all the estimates. However, this does not change the general structure of the estimate \eqref{gradient-estimate-dphi}.
\item Our calculation shows that the magnitude of \(A(d\phi,d\phi)\) clearly has the strongest influence on the estimate on \(|d\phi|\).
\item For Dirac-harmonic maps gradient estimates have been established in \cite{MR3144359}, the authors used 
a Kato-	Yau inequality to obtain the optimal constants in their estimates.
However, this does not seem to help much here since we are considering a more complicated system as in \cite{MR3144359}.
\end{enumerate}
\end{Bem}

\bibliographystyle{plain}
\bibliography{mybib}
\end{document}